\documentclass{amsproc}
\usepackage{verbatim}
\usepackage{url}
\usepackage{amscd,amssymb,bm}
\usepackage{graphics}

\newtheorem{Thm}{Theorem}[section]
\newtheorem{Lem}[Thm]{Lemma}
\newtheorem*{Thmt}{``Theorem''}
\newtheorem{Prop}[Thm]{Proposition}
\newtheorem{Cor}[Thm]{Corollary}
\newtheorem{Conj}[Thm]{Conjecture}
\newtheorem*{Thm*}{Theorem}
\newtheorem*{Lem*}{Lemma}

\theoremstyle{definition}

\newtheorem{Def}[Thm]{Definition}

\theoremstyle{remark}
\newtheorem{Rem}[Thm]{Remark}
\newtheorem*{Ack}{Acknowledgment}

\newtheorem*{Rem*}{Remark}

\numberwithin{equation}{section}

\newcommand{\ndash}{\nobreakdash-\hspace{0pt}}
\newcommand{\Ndash}{\nobreakdash--}

\newcommand\und{\underline}

\newcommand{\frg}{{\mathfrak{g}}}

\newcommand{\frh}{{\mathfrak{h}}}

\newcommand{\calB}{\mathcal{B}}

\newcommand{\calC}{\mathcal{C}}

\newcommand{\calO}{\mathcal{O}}
\newcommand{\calM}{\mathcal{M}}
\newcommand{\calV}{\mathcal{V}}
\newcommand{\calD}{\mathcal{D}}

\newcommand{\calVe}{\mathcal{V}_\epsilon}
\newcommand{\calDe}{\mathcal{D}_\epsilon}

\newcommand{\ii}{{\mathrm{i}}}
\newcommand{\dd}{{\mathrm{d}}} 
\newcommand{\ee}{{\mathrm{e}}}

\newcommand{\sfA}{{\mathsf{A}}}
\newcommand{\sfB}{{\mathsf{B}}}

\newcommand{\sfM}{{\mathsf{M}}}

\newcommand{\sfC}{{\mathsf{C}}}
\newcommand{\sfI}{{\mathsf{I}}}

\newcommand{\sfX}{{\mathsf{X}}}

\newcommand{\sfeta}{\boldsymbol{\eta}}
\newcommand{\sfxi}{\boldsymbol{\xi}}

\newcommand{\bbR}{{\mathbb{R}}}

\newcommand{\bbZ}{{\mathbb{Z}}}

\newcommand{\de}{\partial}

\newcommand{\vev}[1]{{\left\langle\;{#1}\;\right\rangle}}

\newcommand{\Map}{\operatorname{Map}}

\newcommand{\Poiss}[2]{\left\{{\,{#1}\,,\,{#2}\,}\right\}}
\newcommand{\Lie}[2]{\left[{\,{#1}\,,\,{#2}\,}\right]}
\newcommand{\sfPoiss}[2]{\bm{\left\{}{\,{#1}\,,\,{#2}\,}\bm{\right\}}}

\newcommand{\sfstar}{{\bm\star}}
\newcommand{\sfLie}[2]{{\bm{\left[}{\,{#1}\,,\,{#2}\,}\bm{\right]}}}

\DeclareMathOperator{\Ima}{Im}
\DeclareMathOperator{\Ker}{Ker}
\DeclareMathOperator{\codim}{codim}
\DeclareMathOperator{\Graph}{Graph}

\newcommand\qq{}
\newcommand\cmp[1]{{\qq Commun.\ Math.\ Phys.\ \bf #1}}

\newcommand\pl[1]{{\qq Phys.\ Lett.\ \bf #1}}

\newcommand\mpl[1]{{\qq Mod.\ Phys.\ Lett.\ \bf #1}}

\newcommand\lmp[1]{{\qq Lett.\ Math.\ Phys.\ \bf #1}}

\newcommand\ijmp[1]{{\qq Int.\ J. Mod.\ Phys.\ \bf #1}}

\newcommand\anp[1]{{\qq Ann.\ Phys.\ \bf #1}}

\newcommand\adm[1]{{\qq Adv.\ Math.\ \bf #1}}

\newcommand\jdg[1]{{\qq J.\ Diff.\ Geom.\ \bf #1}}

\newcommand\dmj[1]{{\qq Duke Math.\ J. \bf #1}}

\newcommand\atmp[1]{{\qq Adv.\ Theor.\ Math.\ Phys.\ \bf #1}} 
\newcommand\jpaa[1]{{\qq J. Pure Appl.\ Algebra \bf #1}}
\newcommand\invm[1]{{\qq Invent.\ Math.\ \bf #1}}
\newcommand\travm[1]{{\qq Travaux math\'ematiques \bf #1}}
\newcommand\jms[1]{{\qq J. Math. Sci. \bf #1}}

\begin{document}

\title{Deformation Quantization and Reduction}

\author{Alberto~S.~Cattaneo}
\address{Institut f\"ur Mathematik, Universit\"at Z\"urich--Irchel,  
Winterthurerstrasse 190, CH-8057 Z\"urich, Switzerland}

\email{alberto.cattaneo@math.unizh.ch}
\thanks{The author acknowledges partial support of SNF Grant \#20-113439. This work has been partially supported
by the European Union through the FP6 Marie Curie RTN ENIGMA (Contract
number MRTN-CT-2004-5652) and by the European Science Foundation through the MISGAM program.}






\begin{abstract}
This note is an overview of the Poisson sigma model (PSM) and its applications in deformation quantization.
Reduction of coisotropic and pre-Poisson submanifolds, their appearance as branes of the PSM,
quantization in terms of $L_\infty$- and $A_\infty$\ndash algebras, and bimodule structures are recalled. 
As an application, an ``almost'' functorial quantization of Poisson maps is presented if no anomalies occur. This leads in principle to a novel approach
for the quantization of Poisson--Lie groups.
\end{abstract}

\maketitle

\section{Introduction}
The Poisson sigma model (PSM) \cite{I,SS} is a two-dimensional topological field theory with target a Poisson manifold $(M,\pi)$.
It is defined by the action functional
\begin{equation}\label{e:PSM}
S=\int_\Sigma \eta\,\dd X +\frac12\pi(X)\eta\,\eta,
\end{equation}
where the pair $(X,\eta)$ is a bundle map $T\Sigma\to T^*M$, and $\Sigma$ is a two-manifold.

The Hamiltonian study of the PSM on $\Sigma=[0,1]\times\bbR$ leads \cite{CF01} to the construction of the symplectic groupoid of $M$ whenever $T^*M$ is
an integrable Lie algebroid. The functional-integral perturbative quantization yields \cite{CF00} Kontsevich's deformation quantization \cite{K}
if $\Sigma$ is the disk and one chooses vanishing boundary conditions for $\eta$.

General boundary conditions (``branes'')
compatible with symmetries and the perturbative expansion were studied in \cite{CF04} and turned out to correspond to
coisotropic submanifolds. In \cite{C,CF05} (see also references therein) general boundary conditions compatible just with symmetries were classified; they
correspond to what we now call pre-Poisson submanifolds (a generalization of the familiar notion of pre-symplectic submanifolds in the symplectic case).

The perturbative functional-integral quantization with branes gives rise to deformation quantization of the reduced spaces \cite{CF04,CF05}.
The many-brane case leads to the construction of bimodules and morphisms between them \cite{CF04} and potentially to a method for quantizing
Poisson maps in an almost functorial way (see Section~\ref{s:MB}).
This problem is also discussed in \cite{B04} where results (and obstructions) for the symplectic case are given.
The linear case (i.e., the case when the Poisson manifold is the dual of a Lie algebra and the coisotropic submanifold is an affine subspace) leads
to many interesting results in Lie theory \cite{CF04,CT}.

A more general approach leads to a cohomological description of coisotropic submanifolds by $L_\infty$\ndash algebras naturally associated to them 
\cite{OP} and to a deformation quantization in terms of (possibly nonflat) $A_\infty$\ndash algebras. A natural way to reinterpret these results
is by associating a graded manifold to the coisotropic submanifold and use a graded version of Kontsevich's $L_\infty$\ndash quasi\ndash isomorphism
\cite{CF07,LS}. This may also be regarded as a duality for the PSM with target a graded manifold \cite{CF07,C06}.

In the absence of the so-called anomaly---i.e., when the $A_\infty$\ndash algebra can be made flat---one can go down to cohomology.

This paper is an overview of all these themes. It contains a more extended discussion, Section~\ref{s:ft}, of deformation
quantization for graded manifolds and an introduction, Section~\ref{s:MB}, to methods for an almost functorial quantization of morphisms whenever the anomaly
is not present.

\subsection*{Plan of the paper}
In Section~\ref{s:red} we recall properties and reduction of coisotropic, presymplectic and pre-Poisson submanifolds
collecting results in \cite{C,CF05,CZ,CZbis}.
In Section~\ref{s:cd} we give a simple derivation of the $L_\infty$\ndash structure 
associated \cite{OP} to a coisotropic submanifold; we also recall the BFV 
(Batalin--Fradkin--Vilkovisky) formalism
\cite{BF,BV} and its relation \cite{S} to the $L_\infty$\ndash structure.
In Section~\ref{s:redPSM} we recall the Hamiltonian description of the PSM and how coisotropic \cite{CF04} and pre-Poisson \cite{C,CF05}
submanifolds show up as boundary conditions (branes) for the PSM\@. Perturbative functional-integral quantization of the PSM \cite{CF00} is recalled in
Section~\ref{s:exp} with results for coisotropic \cite{CF04} and pre-Poisson \cite{CF05} branes.
Section~\ref{s:sPSM} deals with the reinterpretation of these results in terms of a duality for the PSM with target a graded manifold \cite{CF07,C06}.
In Section~\ref{s:ft} we recall the Formality Theorem for graded manifolds \cite{K,CF07} and its applications in deformation quantization \cite{CF07,LS};
the potential anomaly\cite{CF04,CF07} is discussed;
subsection~\ref{s:qc} contains a description
of the induced properties in cohomology and is original; subsection~\ref{s:Psub} contains a
new description of methods for quantizing the inclusion map of a Poisson submanifold; in particular, it contains a proof of a conjecture by
\cite{CR} that the deformation quantization of a Poisson submanifold determined by central constraints may be obtained by quotienting Kontsevich's
deformation quantization by the ideal generated by the image of the constraints under the
Duflo--Kirillov--Kontsevich map.
Finally, Section~\ref{s:MB} recalls the many-brane case \cite{CF04} and, in the absence of anomalies,
presents a novel method for the quantization of Poisson maps compatible with compositions; the application to the quantization
of Poisson--Lie groups is briefly described.



\begin{Ack}
I thank G.~Felder, F.~Sch\"atz, J.~Stasheff and M.~Zambon for useful discussions and comments.
I thank Hans-Christian Herbig for pointing out reference \cite{B00} (which also led to an improvement
of \cite{S}).
\end{Ack}

\section{Reduction}\label{s:red}
We recall reduction in the context of symplectic and Poisson manifolds. We will essentially follow \cite{CZ} (see also \cite{CZbis}).
Let us first recall the linear case.

A \textsf{symplectic space} is a (possibly infinite dimensional) linear space $V$ equipped with a nondegenerate skew-symmetric bilinear form $\omega$, called
a \textsf{symplectic form}.
By nondegenerate we mean here that $\omega(v,w)=0\ \forall v\in V\Rightarrow w=0$ (this is sometimes called a weak symplectic form). We denote by $\omega^\sharp$
the induced linear map $V\to V^*$: $v\mapsto \omega(v,\ )$. With this notation, $\omega$ is nondegenerate if{f} $\omega^\sharp$ is injective.

The restriction of $\omega$ to a subspace $W$ is in general degenerate, the kernel of $(\omega|_W)^\sharp$ being $W^\perp\cap W$,
where
\[
W^\perp :=\{v\in V : \omega(v,w)=0\ \forall w\in W\}
\]
is the \textsf{symplectic orthogonal} to $W$.
Thus, the \textsf{reduced space} $\und W:=W/(W\cap W^\perp)$ is automatically endowed with a symplectic form. 

An extreme case is when $W\cap W^\perp=0$; i.e., when $(\omega|_W)^\sharp$ is nondegenerate. In this case, $W$ is called a \textsf{symplectic subspace}
and $\und W=W$. In finite dimensions, one has the additional properties that $W^\perp$ is also symplectic and that
\begin{equation}\label{e:cosym}
W\oplus W^\perp = V.
\end{equation}
That is, in finite dimensions a symplectic subspace may be equivalently described as a subspace which admits a symplectic complement. For  this reason,
it may as well be called a \textsf{cosymplectic subspace}.

The other nontrivial extreme case is when $W^\perp\subset W$. In this case $W$ is called a \textsf{coisotropic subspace} and
$\und W=W/W^\perp$. The coisotropic case is in some sense the only case one has to consider, as we have the following simple
\begin{Lem}\label{Lsubsp}
Let $W$ be a subspace of a finite dimensional symplectic space $V$. 
Then it is always possible to find a symplectic subspace $W'$ of $V$ which contains $W$ as a coisotropic
subspace.
\end{Lem}
The Lemma is proved by adding $W$ what is missing to make it symplectic. Any complement to $W+W^\perp$ will do (see Lemma 3.1 in \cite{CZ} for details).

A \textsf{symplectic manifold} is a smooth manifold endowed with a closed, nondegenerate $2$\ndash form. Here nondegenerate means that the bundle map
$\omega^\sharp\colon TM\to T^*M$ induced by the $2$\ndash form $\omega$ on $M$ is injective. The restriction of $\omega$ to a submanifold $C$ is in general
degenerate. The kernel of $(\omega|_C)^\sharp$ is the \textsf{characteristic distribution} $T^\perp C\cap TC$ which in general is not smooth. Here $T^\perp C$ 
is the bundle of symplectic
orthogonal spaces to $TC$ in the restriction
$T_CM$ of $TM$ to $C$. 
The submanifold $N$ is called \textsf{presymplectic} if its characteristic distribution is smooth 
(i.e., $T^\perp C\cap TC$ is a subbundle of $TC$); by the closedness of $\omega$ it is automatically integrable. 
The corresponding leaf space $\und C$ is called the \textsf{reduced space}.
If it is a manifold, it is endowed with a unique symplectic form whose pullback to $C$ is equal to the restriction of $\omega$. The extreme nontrivial examples
of presymplectic submanifolds are the \textsf{symplectic} and \textsf{coisotropic} submanifolds. 
A submanifold $C$ is called
\textsf{symplectic} (\textsf{coisotropic}) if $T_xC$ is symplectic (coisotropic) in $T_xM$ for every $x\in C$. 
The reduction of a symplectic submanifold is the manifold itself. If $C$ is a  submanifold of a finite dimensional manifold $M$, by Lemma~\ref{Lsubsp}
we may find a symplectic subspace of $T_xM$ containing $T_xC$ as a coisotropic subspace for every $x\in C$. These subspaces may be chosen to be glued
together smoothly if the submanifold is presymplectic. Namely, we have the following
\begin{Prop}[\cite{CZ}]\label{t:presympl}
Let $C$ be a presymplectic submanifold of a symplectic manifold $M$. Then it is always possible to find a symplectic submanifold which contains $C$
as a coisotropic submanifold.
\end{Prop}
Moreover, one can show that these symplectic extensions are neighborhood equivalent. Namely, given two such extensions $C'$ and $C''$, there exists a tubular
neighborhood $U$ of $C$ such that $U\cap C'$ and $U\cap C''$ are related by a symplectomorphisms (i.e., a diffeomorphism compatible with the symplectic form)
of $U$ which fixes $C$.

We now move to the Poisson case. We will only work in finite dimensions.
A \textsf{Poisson space} is a linear space $V$ endowed with a bivector $\pi$ (i.e., an element of $\Lambda^2V$). We denote by $\pi^\sharp$
the induced linear map $V^*\to V$, $\alpha\mapsto\pi(\alpha,\ )$. A finite dimensional symplectic space $(V,\omega)$ is automatically Poisson with
$\pi^\sharp=(\omega^\sharp)^{-1}$. The Poisson generalization of the notion of symplectic orthogonal of a subspace $W$ is the image
under $\pi^\sharp$ of its annihilator $W^0:=\{\alpha\in V^*:\alpha(w)=0\ \forall w\in W\}$. In the symplectic case, $\pi^\sharp(W^0)=W^\perp$.
Condition \eqref{e:cosym} may be replaced by $W\oplus\pi^\sharp(W^0)=V$; a subspace $W$ satisfying it is called 
\textsf{cosymplectic}; it may equivalently be characterized by the condition that
the projection of $\pi$ to $\Lambda^2(V/W)$ is symplectic. A \textsf{coisotropic subspace} is analogously defined as a subspace $W$ with $\pi^\sharp(W^0)\subset W$.

A Poisson manifold is a manifold $M$ endowed with a bivector field $\pi$ (i.e., a section of $\Lambda^2TM$) such that the bracket $\Poiss fg:=\pi(\dd f,\dd g)$,
$f,g\in C^\infty(M)$ satisfies the Jacobi identity. Equivalently, $\Lie\pi\pi=0$ with the Schouten--Nijenhuis bracket. This also amounts to saying that
$(C^\infty(M),.,\Poiss{\ }{\ })$ is a Poisson algebra. We denote by $\pi^\sharp\colon T^*M\to TM$ the corresponding bundle map. A finite dimensional symplectic manifold
$(M,\omega)$ is Poisson with $\pi^\sharp={\omega^\sharp}^{-1}$. A submanifold $C$ is called cosymplectic (coisotropic) if $T_xC$ is so in $T_xM$ $\forall x\in C$.
We will denote by $N^*C$---the \textsf{conormal bundle}---the vector bundle over $C$ with fiber at $x$ given by $N^*_xC:=(T_xC)^0\subset T^*_xM$.
Observe that $\pi^\sharp(N^*C)$ is a (singular) distribution on $C$ if $C$ is coisotropic. Invariant functions form naturally a Poisson algebra. If the leaf space
$\und C$ is a manifold, it is naturally a Poisson manifold. A cosymplectic manifold $C$ is also automatically a Poisson manifold (not a Poisson submanifold though).
One way to see this is to think of a Poisson manifold as a manifold foliated by symplectic leaves. A cosymplectic submanifold intersects symplectic leaves
cleanly, and each intersection is a symplectic submanifold of the symplectic leaf; thus, a cosymplectic submanifold
is foliated in symplectic leaves, which makes it into a Poisson manifold. The induced Poisson bracket may also be obtained by Dirac's procedure.



We now wish to describe a Poisson generalization of the notion of presymplectic submanifold leading to a generalization of Proposition~\ref{t:presympl}.
Namely, we call a submanifold $C$ \textsf{pre-Poisson} if $\pi^\sharp(N^*C)+TC$ has constant rank along $C$. (In \cite{CF05} and references therein,
this was called a ``submanifold with strong regularity constraints''.) In the symplectic case, this is equivalent to $C$ being presymplectic.

An equivalent definition of a pre-Poisson submanifold $C$ amounts to asking that the bundle map $\phi\colon N^*C\to NC$ obtained by composing
the injection $N^*C\to T^*_CM$ with $\pi^\sharp$ and finally with the projection $T_CM\to NC$ should have constant rank. As special cases we recognize
coisotropic submanifolds ($\phi=0$) and cosymplectic submanifolds ($\phi$ surjective). We have the following
\begin{Prop}[\cite{CF05,CZ}]\label{p:CF}
Let $C$ be a pre-Poisson submanifold of a Poisson manifold $M$. Then it is always possible to find a cosymplectic submanifold $M'$ which contains $C$
as a coisotropic submanifold. Moreover, this extension $M'$ is unique up to neighborhood equivalence.
\end{Prop}

In the following we will also need a description in terms of Lie algebroids. Recall that a \textsf{Lie algebroid} is a vector bundle $E$ over a smooth
manifold $M$ endowed with a bundle map $\rho\colon E\to TM$ (the anchor map) and a Lie algebra structure on $\Gamma(E)$ with the property
$\Lie a{fb}=f\Lie ab + \rho(a)f\,b$ for every $a,b\in\Gamma(E)$, $f\in C^\infty(M)$. The tangent bundle itself is a Lie algebroid.
The cotangent bundle of Poisson manifold is also a Lie algebroid
with $\rho=\pi^\sharp$ and $\Lie{\dd f}{\dd g}:=\dd\Poiss fg$. If $C$ is a pre-Poisson submanifold, then $AC:=N^*C\cap{\pi^\sharp}^{-1}TC$ is a Lie 
subalgebroid \cite{CF05}. Its anchor is defined just by restriction; as for the bracket, one has to extend sections of $AC$ to sections of $T^*M$ in 
a neighborhood of $C$, use the bracket on $T^*M$ and finally restrict; it is not difficult to check that the result is independent of the extension.
Observe that for $C$  coisotropic, one simply has $AC=N^*C$.

\section{Cohomological descriptions}\label{s:cd}
The complex $\Gamma(\Lambda TM)$ of multivector fields  on a smooth manifold $M$ is naturally endowed with a Lie bracket of degree $-1$, the 
\textsf{Schouten--Nijenhuis bracket}. Let $\calV(M):=\Gamma(\Lambda TM)[1]$ be the corresponding graded Lie algebra (GLA).
A Poisson bivector field is now the same as an element $\pi$ of $\calV(M)$ of degree one which satisfies $\Lie\pi\pi=0$.
In general, a self-commuting element of degree $1$ in a GLA is called a Maurer--Cartan (MC) element. Observe that $\Lie\pi{\ }$ is a differential; actually,
it is the Lie algebroid differential corresponding to the Lie algebroid structure on $T^*M$.

Let $C$ be a submanifold and $NC$ its normal bundle, canonically defined as the quotient by $TC$ of the restriction $T_CM$ of $TM$ to $C$. It is also the dual
bundle of the conormal bundle $N^*C$ introduced above.
Let $A$ be the graded commutative algebra $\Gamma(\Lambda NC)$.
Denote by $P\colon\Gamma(\Lambda TM)\to A$ 
the composition of the restriction to $C$ with the projection $T_CM\to NC$.
It is not difficult to check that $\Ker P$ is a Lie subalgebra. Moreover, $C$ is coisotropic if{f} $\pi\in\Ker P$. In this case, there is a 
well-defined differential $\delta$ on $A$ acting on $X\in A$ by
\[
\delta X:=P\Lie\pi{\Tilde X},\quad \Tilde X\in P^{-1}(X).
\]
This is the differential corresponding to the Lie algebroid structure on $N^*C$. Observe that $H^0_\delta$ is the algebra of invariant functions on $C$.

The projection $P$ admits a section if $M$ is a vector bundle over $C$. In particular, this is true if $M$
is the normal bundle of $C$.
The section $i$ for $P\colon\Gamma(\Lambda T(NC))\to A$ is given as follows: If $f$ is a function on $C$, one defines $if$ as the pullback
of $f$ by the projection $NC\to C$. If $X$ is a section of $NC$, one defines $iX$ as the unique vertical vector field $\Tilde X$ on $NC$ which is constant along
the fibers and such that $P\Tilde X = X$. Finally, $i$ is extended to sections of $\Lambda NC$ so that it defines an algebra morphism. It turns out that
$i(\Gamma(NC))$ is an abelian subalgebra.

Let us now choose an embedding of $NC$ into our Poisson manifold $(M,\pi)$. By restriction $\pi$ defines a Poisson structure on $NC$, so a MC element in
$\calV(NC)$. Following Voronov \cite{V}, we define the derived brackets
\begin{equation}\label{e:lambda}
\begin{split}
\lambda_k&\colon A^{\otimes k}\to A,\\
\lambda_k(a_1,\dots,a_k)&:=
P\left(\Lie{\Lie{\dots\Lie{\Lie \pi{ia_1}}{ia_2}}\dots}{i a_k}\right),
\end{split}
\end{equation}
and $\lambda_0:=P(\pi)\in A$. Observe that $\lambda_0=0$ if{f} $C$ is coisotropic; in this case, $\delta=\lambda_1$.
In general Voronov proved the following
\begin{Thm}\label{t:Vor}
Let $\frg$ be a GLA, $\frh$ an abelian subalgebra, $i$ the inclusion map. Let $P$ be a projection to $\frh$ such that $P\circ i=\mathit{id}$
and that $\Ker P$ is a Lie subalgebra. Then the derived brackets of every MC element of $\frg$ define an $L_\infty$\ndash structure on $\frh$.
\end{Thm}
This means that the operations $\lambda_k$ satisfy certain quadratic relations which in particular for $\lambda_0=0$ imply that $\lambda_1$ is a differential
for $\lambda_2$ and that $\lambda_2$ satisfies the Jacobi identity up to $\lambda_1$\ndash homotopy. So the $\lambda_1$\ndash cohomology inherits the structure
of a GLA.

In our case $A$ is also a graded commutative algebra and the multibrackets are multiderivations. In this case we say that we have a $P_\infty$\ndash algebra
($P$ for Poisson). In particular, when $C$ is coisotropic, $H_{\lambda_1}=H_\delta$ is a graded Poisson algebra. The Poisson structure in degree zero---i.e.,
on invariants functions---is the same as the one described in the previous Section.

Observe that the $P_\infty$\ndash structure depends on a choice of embedding $NC\hookrightarrow M$. It is possible to show that different choices lead to
$L_\infty$\ndash isomorphic algebras. We will return on this in a forthcoming paper \cite{CS}.

The $P_\infty$\ndash structure appeared first in \cite{OP} as a tool to describe deformations of a coisotropic submanifold. 
In \cite{CF07} it was rediscovered as the semiclassical limit of the quantization of coisotropic submanifolds. We will return on quantization in Section~\ref{s:ft}.


\subsection{The BFV method}
The Batalin--Fradkin--Vilkovisky (BFV) method is an older method to describe coisotropic submanifold in terms of a differential graded Poisson 
algebra (DGPA) \cite{BF,BV}. 
Its advantage is that DGPAs are more manageable than $P_\infty$\ndash algebras. The disadvantage is that it works only under certain
assumptions. We will essentially follow Stasheff's presentation \cite{S97}.

Suppose first that the coisotropic submanifold $C$ is given by global constraints, or, equivalently, that the normal bundle of $C$ is trivial.
Let $\{y^\mu\}_{\mu=1,\dots,n}$, $n=\codim C$, be a basis of constraints (equivalently, a basis of sections of $NC$ or a basis of transverse coordinates).
Let $\calB$ be the free graded commutative algebra generated by
a set $\{b^\mu\}$ of
odd variables  of degree $1$ (the ``ghost momenta'').
On $C^\infty(NC)\otimes\calB$ one defines the Koszul differential
$\delta_0:=y^\mu\frac\de{\de b^\mu}$. The $\delta_0$ cohomology is then concentrated in degree $0$ and is $C^\infty(C)$. Let us also introduce a set 
$\{c_\mu\}$ of odd variables
of degree $-1$ (the ``ghosts'') and the free graded commutative algebra $\calC$ they generate. On $C^\infty(NC)\otimes\calB\otimes\calC$ one has a unique  Poisson
structure with $\Poiss{b^\mu}{c_\nu}=\Poiss{c_\nu}{b^\mu}=\delta^\mu_\nu$. Then $\delta_0$ is a Hamiltonian vector field with Hamiltonian function
$\Omega_0=y^\mu c_\mu$. On this enlarged algebra the $\delta_0$\ndash cohomology is $C^\infty(NC)\otimes\calC\simeq\Gamma(\Lambda NC)$.

By choosing an embedding $NC\hookrightarrow M$, we get a Poisson structure on $C^\infty(NC)$ which we can extend to $C^\infty(NC)\otimes\calB\otimes\calC$.
We now consider the sum of the two Poisson structure. Let $F_0:=\frac12\Poiss{\Omega_0}{\Omega_0}=\frac12\Poiss{y^\mu}{y^\nu}c_\mu c_\nu$.
Observe that $\{F_0,\delta_0,\Poiss{\ }{ }\}$ defines a nonflat $L_\infty$\ndash structure\footnote{Originally
$L_\infty$\ndash algebras were always assumed to miss the 0th term.
When $L_\infty$\ndash algebras with the 0th term turned out to be interesting,
various terminologies were introduced. An $L_\infty$\ndash algebra 
without the 0th term is called nowadays flat (resp.\ strict), while one where
the 0th term is there is called curved (resp.\ weak). If no assumption on the
0th term is made (i.e., it may or may not vanish), we simply speak of an
$L_\infty$\ndash algebra.} 
on $C^\infty(NC)\otimes\calB\otimes\calC$.
Moreover, $F_0$ is $\delta_0$\ndash exact if{f} $C$ is coisotropic. In this case, one can use cohomological perturbation theory to kill $F_0$ as follows.

Let $h:=b^\mu\frac\de{\de y^\mu}$. Then $\Lie{\delta_0}h=E:=b^\mu\frac\de{\de b^\mu}+y^\mu\frac\de{\de y^\mu}$. One then has a homotopy $s$ for $\delta_0$
(i.e. $s\delta_0+\delta_0s=\mathit{id}-\mathit{pr}$, with $\mathit{pr}$ the projection to a subspace isomorphic to the $\delta_0$\ndash cohomology)
defined by $s\alpha=h\alpha/|\alpha|$ if $E\alpha=|\alpha|\alpha$ and $s\alpha=0$ if $E\alpha=0$. Let us then define $\Omega:=\sum_{i=0}^\infty\Omega_i$ by induction 
as follows: Let $R_k:=\sum_{i=0}^k\Omega_i$. Then define $\Omega_{k+1}=-\frac12\,s\left(\Poiss{R_k}{R_k}\right)$.
If $C$ is coisotropic, the following hold:
\begin{enumerate}
\item $\Poiss\Omega\Omega=0$, so  $\{D,\Poiss{\ }{ }\}$ defines a DGPA structure on $C^\infty(NC)\otimes\calB\otimes\calC$, where $D:=\Poiss\Omega{\ }=\delta_0+\cdots$.
\item The $D$\ndash cohomology $H_D$ is isomorphic to the Lie algebroid cohomology $H_\delta$ of $N^*C$.
\item This way $H_\delta$ gets the structure of a GPA.
\end{enumerate}


It is natural to ask whether the GPA structure on $H_\delta$ is the same as the one induced before. The answer is affirmative. Actually, there is a stronger result:
\begin{Thm}[Sch\"atz \cite{S}]
$(C^\infty(NC)\otimes\calB\otimes\calC, D,\Poiss{\ }{\ })$ is $L_\infty$\ndash quasi\ndash isomorphic to $(A,\{\lambda_k\})$.
\end{Thm}

Observe that the $L_\infty$\ndash structure $(A,\{\lambda_k\})$ can be defined for every submanifold, while the BFV formalism as presented above requires
that the normal bundle should be trivial, or at least flat. A more general version of the BFV formalism allows for linearly dependent global constraints
and uses the Koszul--Tate resolution. There is however a generalization by Bordemann and Herbig \cite{B00}
of the construction presented above which needs no assumption on the normal bundle and just requires the choice of a connection. The above Theorem
holds also in the general case \cite{S}.




\section{The reduced space of the Poisson sigma model}\label{s:redPSM}
Given a finite dimensional manifold $M$, we denote by $PM$ the Banach manifold of $C^1$\ndash paths in $M$; viz., $PM:=C^1(I,M)$, $I=[0,1]$. 
We denote by $T^*PM$ the vector bundle over $PM$ with fiber at $X\in PM$ the space of $C^0$ sections of $T^*I\otimes X^*T^*M$.
Using integration over $I$ and the canonical pairing between $TM$ and $T^*M$, one can endow $T^*PM$ with a symplectic structure.
Let now $\pi$ be Poisson bivector field on $M$. We define
\[
\calC(M,\pi):=\{(X,\zeta)\in T^*PM : \dd X +\pi^\sharp(X)\zeta = 0\}.
\]

These infinite-dimensional  manifolds naturally appear when considering the PSM \eqref{e:PSM}
on $\Sigma=I\times\bbR$. Namely, the map $X\colon\Sigma\to M$
may be regarded as a path in $PM$. On the other hand, introducing the coordinate $t$ on $\bbR$, we may write $\eta=\zeta+\lambda\dd t$,
and reinterpret the pair $(X,\zeta)$ as a path in $T^*PM$. The bilinear term in $X$ and $\zeta$ in the
action functional \eqref{e:PSM} yields the canonical symplectic structure on $T^*PM$, while $\lambda$ appears linearly as a Lagrange multiplier
leading to the constraints that define $\calC(M,\pi)$.

\begin{Thm}[\cite{CF01}]
$\calC(M,\pi)$ is a coisotropic submanifold of $T^*PM$. Its reduced space $\und\calC(M,\pi)$ is
the source simply connected symplectic groupoid of $M$, whenever $M$ is integrable (i.e., $T^*M$ is an integrable Lie algebroid). 
\end{Thm}
We now want to discuss boundary conditions. Given two submanifolds $C_0$ and $C_1$ of $M$, let
\[
\calC(M,\pi;C_0,C_1):=\{(X,\zeta)\in\calC(M,\pi) : X(0)\in C_0,\ X(1)\in C_1\}.
\]
Observe that $\calC(M,\pi;C_0,C_1)$ has natural maps to $C_0$ and $C_1$ (evaluation of $X$ at $0$ and $1$). We call points in the image in $C_0\times C_1$
connectable.

Its symplectic orthogonal bundle may be explicitly computed. For simplicity we choose local coordinates on $M$. (The correct, but more cumbersome, description
involves choosing a torsion-free connection for $TM$.) 
One gets \cite{C}
\begin{multline}
T^\perp_{(X,\zeta)}\calC(M,\pi;C_0,C_1)=
\{(-\pi^\sharp(X)\beta,\dd\beta_i + \de_i\pi^{jk}(X)\zeta_j\beta_k) :\\
 \beta\in\Gamma(X^*T^*M),\ \beta(0)\in N^*_{X(i)}C_i,\ i=0,1\}.
\end{multline}
Intersection with $T\calC(M,\pi;C_0,C_1)$ just forces the boundary conditions for $\beta$ to be such that the $X$-variations are tangent to the
submanifolds: viz., $\pi^\sharp(X(i))\beta(i)\in T_{X(i)}C_i$, $i=0,1$. Thus, $\calC(M,\pi;C_0,C_1)$ is presymplectic if{f} the pre-Poisson condition
is satisfied at all connectable points \cite{CF05,C}, and it is coisotropic if{f} the coisotropicity condition is satisfied at all connectable points \cite{CF04}.
The sections of $\delta X\oplus\delta\zeta$ of the characteristic distribution
are parametrized by a section $\beta$ of $X^*T^*M$ with $\beta(i)\in A_{X(i)}C_i$, $i=0,1$:\footnote{Actually, this is the image of the anchor of an infinite-dimensional
Lie algebroid over $\calC(M,\pi;C,M)$. Explicitly this is described in \cite{BC} and in \cite{BCZ}.}
\begin{align}\label{e:deltaHam}
\delta X&=-\pi^\sharp(X)\beta,\\
\delta\zeta_i&=\dd\beta_i + \de_i\pi^{jk}(X)\zeta_j\beta_k,
\end{align}
where for simplicity we have written the second equation using local coordinates on $M$. (A more invariant, but more cumbersome, possibility
is to write it upon choosing a torsion-free connection for $TM$ and observing that the distribution does not depend on this choice.)

To remove the condition on connectable points, one can e.g.\ take the second submanifold to be the whole manifold $M$. So we have
\begin{Thm}~
\begin{enumerate}
\item $\calC(M,\pi;C,M)$ is coisotropic in $T^*PM$ if{f} $C$ is coisotropic in $M$ \cite{CF04}.
\item $\calC(M,\pi;C,M)$ is presymplectic in $T^*PM$ if{f} $C$ is pre-Poisson in $M$ \cite{CF05,C}.
\end{enumerate}
\end{Thm}
Point 2) says in particular that pre-Poisson submanifolds are the most general boundary conditions for the Poisson sigma model compatible with symmetries.
A submanifold chosen as a boundary condition is usually called a \textsf{brane}.
Point 1) is important as it describes boundary conditions compatible with symmetries \emph{and}\/ perturbation theory around $\pi=0$. 
Namely, consider the family of Poisson structures $\pi_\epsilon:=\epsilon \pi$, $\epsilon\ge0$. 
Let $C$ be a pre-Poisson submanifold for $\pi$. Then it is pre-Poisson for $\pi_\epsilon$, $\forall\epsilon\ge0$. In particular, it is coisotropic for $\epsilon=0$.
One may indeed check that the codimension of the characteristic distribution of $\calC(M,\pi_\epsilon;C,M)$ is the same for all $\epsilon>0$ but it may change
for $\epsilon=0$. It stays the same if{f} $C$ is coisotropic for $\pi$ (and hence coisotropic for $\pi_\epsilon$ for all $\epsilon$).




\section{Expectation values in perturbation theory}\label{s:exp}
We now consider the perturbative functional-integral quantization of the Poisson sigma model. Namely, we consider integrals of the form
\begin{equation}\label{e:vev}
\vev\calO:=\int \ee^{\frac\ii\hbar S}\,\calO,
\end{equation}
where $\calO$ is a function of the fields and the integration is over all fields $(X,\eta)$, $X\in\Map(\Sigma,M)$, $\eta\in\Gamma(T\Sigma\otimes X^*T^*M)$.
If $\Sigma$ has a boundary, the discussion in Section~\ref{s:redPSM} shows that we have to choose a pre-Poisson submanifold $C$ of $M$ (a brane)
and impose
the conditions $X(\de\Sigma)\subset C$, $\iota_\de^*\eta\in\Gamma(T^\Sigma\otimes X^*AC)$, where $\iota_\de$ is the inclusion map $\de\Sigma\hookrightarrow\Sigma$.

Also observe that the Poisson sigma model has symmetries 
(see~\eqref{e:deltaHam}). This makes the integrand $\ee^{\frac\ii\hbar S}$ invariant under the Lagrangian extension of the symmetries which has the form:
\begin{align}
\delta X&=-\pi^\sharp(X)\beta,\\
\delta\eta_i&=\dd\beta_i + \de_i\pi^{jk}(X)\eta_j\beta_k,
\end{align}
with $\beta\in\Gamma(X^*T^*M)$, $\iota_\de^*\beta\in\Gamma(X^*AC)$.
As a consequence the integral \eqref{e:vev} cannot possibly converge (it would not even if we were in finite dimensions). The trick is to select representatives
of the fields modulo symmetries in a consistent way. This is done using the BV formalism \cite{BV} which in particular guarantees that the result
is independent of the choice involved (``gauge fixing''). Its application to the PSM is described
in \cite{CF00}. An important issue is that an expectation value as in \eqref{e:vev} is also gauge-fixing independent if $\calO$ is invariant under symmetries.
We will consider only boundary observables. Namely, let $f$ be a function on $C$. To it we associate the observable $\calO_{f,u}(X,\eta):=f(X(u))$, where
$u$ is some point of $\de\Sigma$.
The observable is invariant if{f} $f$ is invariant constant along
the characteristic distribution of $C$. It will be interesting for us to take two such functions $f$ and $g$
and consider $\calO=\calO_{f,u}\calO_{g,v}$.

We compute integrals as in \eqref{e:vev} in perturbation theory, i.e., by expanding $S$ around its critical points and then computing the
integral formally in terms of the momenta of the Gaussian distribution associated to the Hessian of $S$ (observe that after gauge fixing this is no longer
degenerate). We will pick only the trivial critical point; viz., the class of solutions $\{X$ constant, $\eta=0\}$. We fix the constant by imposing $X$
to be equal to some $x\in C$ at some boundary point which we denote by $\infty$.
Finally, we select $\Sigma$ to be a disk\footnote{Higher genera would be interesting to consider, but this has not been done so far, the main difficulty
being the appearance of nontrivial cohomological classes in the space of solutions.} and compute
\begin{equation}\label{e:starvev}
f\star g(x;u,v):=\int_{X(\infty)=x}\ee^{\frac\ii\hbar S}\,\calO_{f,u}\calO_{g,v},
\end{equation}
where the points $u,v,\infty$ are cyclically ordered. It turns out that such an integral may be explicitly computed in perturbation theory.
The results are as follows.
\begin{figure}[ht] 
\begin{center}
\resizebox{3 cm}{!}
{\includegraphics{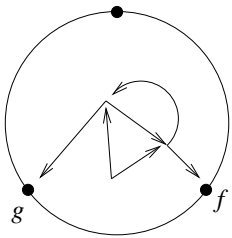}}
\caption{An example of a Kontsevich graph}  
\end{center}
\end{figure}

The case $C=M$ has been studied in \cite{CF00}.
The result turns out to be independent of $u$ and $v$ and to define a star product \cite{BFFLS} on $M$: viz., 
an associative product $\star$ on $C^\infty(M)[[\epsilon]]$
(in the above computation $\epsilon=\ii\hbar/2$) deforming the pointwise product in the direction of the Poisson bracket:
\begin{equation}\label{e:star}
f\star g= fg + \epsilon\Poiss fg + \sum_{n=2}^\infty B_n(f,g),
\end{equation}
where the $B_n$s are bidifferential operator with the property $B_n(1,\bullet)=B_n(\bullet,1)=0$. Moreover, this star product is exactly the one
defined by Kontsevich in \cite{K}.

The case $C$ pre-Poisson has been considered in \cite{CF05}. It turns out that the PSM on $M$ with boundary conditions on $C$ is equivalent to
the PSM on any $M'$ as in Prop.~\ref{p:CF} with the same boundary conditions. As a result, it is enough to consider coisotropic 
submanifolds.

The case $C$ coisotropic has been considered in \cite{CF04}. A few interesting, but unpleasant phenomena appear. First, it might happen that
the result of \eqref{e:starvev} depends on $u$ and $v$. Second, the result regarded as a function on $C^\infty(C)$ may not be invariant.
Third, associativity may not be guaranteed. Under suitable conditions, one may ensure that everything works. We return on this in  subsection~\ref{s:dqcs},
where we show that the vanishing of the first and second Lie algebroid cohomologies is a sufficient condition.

\section{Super PSM and duality}\label{s:sPSM}  
To deal with symmetries in the PSM consistently, one has to resort to the BV formalism. In the case at hand the recipe (following from the AKSZ
formalism \cite{AKSZ} adapted to the PSM with boundary \cite{CF01b}) consists just of formally replacing the fields $(X,\eta)$ in the action by ``superfields'' $(\sfX,\sfxi)$.
To simplify the discussion, let us choose a coordinate chart on $M$. Then we may regard $X$ as an element of 
$\Omega^0(\Sigma)\otimes\bbR^m$, $m=\dim M$,
and $\eta$ as an element of $\Omega^1(\Sigma)\otimes(\bbR^m)^*$. The action in local coordinate reads
\[
S = \int_\Sigma \eta_i\dd X^i + \frac12 \pi^{ij}(X)\eta_i\eta_j,
\]
where we use Einstein's convention on repeated indices. 

The superfields are then $\sfX\in V:=\Omega(\Sigma)\otimes\bbR^m$,
$\sfeta\in W:=\Omega(\Sigma)\otimes(\bbR^m)^*$ and one just puts them in the action instead of the classical fields. The integral selects by definition
the two-form 
component.\footnote{An invariant description requires introducing graded manifolds and regarding the space of superfields as the graded infinite-dimensional manifold $\Map(T[1]\Sigma,T^*[1]M)$.}
The superaction is a function on $V\oplus W$ regarded as a graded vector space.
Observe that for consistency, we have to consider
elements of $V$ as even and those of $W$ as odd. 
It is also useful to introduce a $\bbZ$\ndash grading and assign $V$ degree $0$ and $W$ degree $1$, while integration over $\Sigma$ is given degree $-2$.

If we only consider the ``kinetic term'' $S_0=\int_\Sigma \sfeta_i\dd \sfX^i$, it is quite natural to observe that there is a duality obtained
by exchanging $\sfeta$ and $\sfX$. The only problem is that we want to think of the zero form in $\sfX$ as a map. However, the zero form in $\sfeta$ has
to be considered as an odd element. The problem is solved if we allow the target $M$ to be a graded manifold itself. The duality then exchanges
$M$ with a dual manifold $\Tilde M$. If we now take into account the term involving $\pi$, we see that it corresponds to a similar term for a multivector
field on $\Tilde M$.

If $\Sigma$ has a boundary, the superfields have to be assigned boundary conditions. We choose a coisotropic submanifold $C$ of $M$ and require
$\sfX$  on the boundary to take values in $C$, while $\sfeta$ has to take values in the pullback of the conormal bundle of $C$.
As we work in perturbation theory, we may actually consider just a formal neighborhood of $C$. Namely, we choose an embedding of the normal
bundle $NC$ of $C$ into $M$ and regard $NC$ as the Poisson manifold. The dual manifold turns out to be the graded manifold $N^*[1]C$.
To $\pi$ there corresponds a multivector field $\tilde\pi$ of total degree $2$. The dual theory can be mathematically understood in terms
of a formality theorem for graded manifolds.

\section{Formality Theorem}\label{s:ft}
Kontsevich's star product is an application of his formality theorem stating that the DGLA of multidifferential operators on a smooth manifold is formal
\cite{K}. His explicit local formula for the $L_\infty$\ndash quasi\ndash isomorphism may also be obtained in terms of expectation values of boundary and bulk
observables in the PSM with zero Poisson structure \cite{CF00}. This construction may be generalized to smooth graded manifolds 
\cite{CF07} (see also \cite{C06} for the PSM version). An application \cite{CF07,LS} is the deformation quantization of coisotropic submanifolds \cite{CF04}.

Let $\calM$ be a smooth graded manifold 
$A=\bigoplus_{k\in\bbZ}A^k$ 
its graded
algebra of functions,
$\calV(\calM)$ its complex of graded multivector fields and $\calD(\calM)$ its
complex of graded multidifferential operators (regarded as a subcomplex of the
Hochschild complex of $A$). 
We consider $\calV(\calM)$ 
as DGLAs by taking the total degree (shifted
by one); viz.:
\begin{align*}
\calV^k(\calM) &:= \{r\text{-multivector fields of homogeneous degree }s
\text{ with }r+s=k+1\},\\
\calD^k(\calM) &:= \{r\text{-multidiff.\ operators 
of homogeneous degree }s
\text{ with }r+s=k+1\}.
\end{align*}
\begin{Thm}[Formality Theorem]
There exists an $L_\infty$\ndash quasi\ndash isomorphism
$U\colon\calV(\calM)\leadsto\calD(\calM)$, hence the DGLA $\calD(\calM)$ is formal.
Moreover, $U$ may be chosen such that
the degree-one component $U_1\colon\calV(\calM)\mapsto\calD(\calM)$ is the
Hochschild--Kostant--Rosenberg (HKR) map.
\end{Thm}
See \cite{K} for the proof in the ordinary case and \cite{CF07} for the case of graded manifolds.

One may also extend $U$ to formal power series in a parameter $\epsilon$
and get an $L_\infty$\ndash quasi\ndash isomorphism
$\calVe(\calM)\leadsto\calDe(\calM)$, where
\[
\calVe(\calM):=\epsilon\calV(\calM)[[\epsilon]],\quad
\calDe(\calM):=\epsilon\calD(\calM)[[\epsilon]].
\]
In particular, since MC elements are mapped to MC elements
by an $L_\infty$\ndash morphism (with no convergence problems in
the setting of formal power series), the Formality Theorem also shows that every MC element of $\calV(\calM)$ gives rise
to a MC element of $\calDe(\calM)$. In the case of ordinary manifolds, the former MC element is the same as a Poisson structure, while the latter
defines a deformation quantization of the former. (A classification theorem also follows by the fact that an $L_\infty$\ndash quasi\ndash isomorphism
actually yields an isomorphism of the sets of MC elements modulo gauge transformations).

The case of graded manifold is important for coisotropic submanifolds. As we recalled in Section~\ref{s:cd}, the algebra $A:=\Gamma(NC)$ naturally appears
in the cohomological description of a coisotropic submanifold $C$ of a Poisson manifold $(M,\pi)$. 
Now one can reinterpret $A$ as the algebra of functions on the graded manifold
$N^*[1]C$ (where $[1]$ denotes a shift by $1$ in the fiber coordinates). Moreover, the restriction of $\pi$ to $NC$
(obtained upon choosing an embedding of $NC$ into $M)$ may be regarded as a MC element on $NC$ but also on $N^*[1]C$. 

One way to see this is via
Roytenberg's Legendre mapping theorem \cite{R} which states the existence of a canonical antisymplectomorphism between
$T^*[n]E$ and $T^*[n](E^*[n])$ for every integer $n$ and for every graded vector bundle $E$. In particular, for $n=1$  and observing that
multivector fields on a graded manifold may be reinterpreted as functions on its tangent bundle shifted by $1$,
one gets an anti-isomorphism
between the GLAs $\calV(E)$ and $\calV(E[1])$; in particular, this yields an isomorphism of the sets of MC elements.
Specializing to $E=NC$ (as a graded vector bundle
concentrated in degree zero) yields the sought for result.

A more direct way is just to observe that each derived bracket $\lambda_k$ in \eqref{e:lambda}
is a multiderivation on the algebra $A$ and so a $k$\ndash vector field on $N^*[1]C$.
Their linear combination is the desired MC element. We now consider a more general situation.

\subsection{Maurer--Cartan elements for graded manifolds}
In general, given a  graded manifold $\calM$ and a MC element $\pi$ in $\calV(\calM)$, 
by Theorem~\ref{t:Vor} (with $\frg=\calV(\calM)$, $\frh=C^\infty(\calM)$, $P$ and $i$ the natural projection and inclusion),
we may construct a $P_\infty$\ndash algebra on $A:=C^\infty(\calM)$
with multibrackets as in \eqref{e:lambda}.
Taking care of degrees, one sees that $\lambda_k$ is a  multiderivation of degree $2-i$, i.e.:
\begin{align*}
\lambda_0 &=P_0 \in A^2,\\
\lambda_1 &\colon  A^j\to A^{j+1},\\
\lambda_2 &\colon  A^{j_1}\otimes A^{j_2} 
\to A^{j_1+j_2},\\
\dots\\
\lambda_i &\colon A^{j_1}\otimes\dots\otimes A^{j_i}
\to A^{j_1+\dots j_i+2-i},\\
\dots
\end{align*}

If $P_0=0$, then $\lambda_0=0$ and we have a flat
$L_\infty$
\ndash algebra. In this case, $\lambda_1$ is a differential and
its cohomology
\[
\sfA^\bullet:=H_{\lambda_1}^\bullet(A)
\]
is a graded Poisson algebra. In the following we will denote by $\sfPoiss{\ }{\ }$
the induced Poisson bracket on $\sfA$.
In particular, $\sfA^0$ is a Poisson algebra (but in general it is not
the algebra of functions on a smooth manifold) whose Poisson bracket we will
simply denote by $\Poiss{\ }{\ }$.


\begin{Def}
A graded manifold $\calM$ with a (flat) $P_\infty$\ndash structure
on $A=C^\infty(\calM)$ is called a \textsf{(flat) $P_\infty$\ndash manifold}.
\end{Def}

We now turn to multidifferential operators.
An element $m$ in $\calDe^1(\calM)$ 
consists of a sequence $m_i$, where
$m_i$ is an
$i$\ndash multidifferential operators of degree $2-i$, i.e.,
\begin{align*}
m_0 &\in \epsilon A^2[[\epsilon]],\\
m_1 &\colon  A^j[[\epsilon]] \to \epsilon A^{j+1}[[\epsilon]],\\
m_2 &\colon  A^{j_1}[[\epsilon]]\otimes A^{j_2}[[\epsilon]] 
\to \epsilon A^{j_1+j_2}[[\epsilon]],\\
\dots\\
m_i &\colon A^{j_1}[[\epsilon]]\otimes\dots\otimes A^{j_i}[[\epsilon]] 
\to \epsilon A^{j_1+\dots j_i+2-i}[[\epsilon]],\\
\dots
\end{align*}
Let us denote by $\chi$ the (extension to $A[[\epsilon]]$) of 
the multiplication on $A$ and define $\mu=\chi+m$; viz.:
\begin{equation}\label{mui}
\mu_i = \begin{cases}
m_i & i\not=2,\\
\chi+m_2 & i=2.
\end{cases}
\end{equation}
The element $m$ is MC, i.e.,
\[
\delta m + \frac12\Lie mm = 0,
\qquad (\delta = \Lie\chi{\ }),
\]
if{f} the operations $\mu_i$ define an $A_\infty$\ndash structure
on $A[[\epsilon]]$. 
If in addition $\mu_0$ vanishes, one says that the $A_\infty$\ndash algebra is flat. In this case,
$H_{\mu_1}^\bullet$ is an associative algebra.

Let $\Tilde\mu$ denote the skew-symmetrization of $\mu$. 
Then $(A[[\epsilon]],\Tilde\mu)$ is an $L_\infty$\ndash algebra.
Since the multiplication $\chi$ is graded commutative, the $\Tilde\mu$
will take values in $\epsilon A[[\epsilon]]$. So, dividing by $\epsilon$ and
working modulo
$\epsilon$, they define operations $\underline{\Tilde\mu_i}$
on $A$ which make it into an $L_{\infty}$\ndash algebra. 
It may be shown, see \cite{CF07},
that the $\underline{\Tilde\mu_i}$s are actually
multiderivations, so we have a $P_\infty$\ndash algebra structure on $A$.

\subsection{Deformation quantization of $P_\infty$\ndash manifolds}
Let $A$ be the algebra of functions on a $P_\infty$\ndash manifold $\calM$.
We denote by $\chi$ the graded commutative product and by $\lambda$ the multibrackets.
\begin{Def}
A deformation quantization of $\calM$ consists of an $A_\infty$\ndash structure
on $(A[[\epsilon]],\mu)$ with $\mu=\chi+m$, such that:
\begin{enumerate}
\item $m$ is of order $\epsilon$;
\item $m$ consists of multidifferential operators;
\item $m$ vanishes when one of its arguments is the unit in $A$;
(i.e., $m_i(a_1,\dots,a_i)=0$ if $a_j=1$ for some $j\in\{1,\dots,i\}$, $i>0$);
\item the induced $P_\infty$\ndash structure $\Tilde\mu$ is equal to $\lambda$.
\end{enumerate}
If $\calM$ is flat, by deformation quantization in the flat sense we mean that in addition the condition
$\mu_0=0$ is fulfilled.
\end{Def}
The formality theorem for graded manifolds then implies the following
\begin{Thm}
Every $P_\infty$\ndash manifold $\calM$ admits a deformation quantization.
\end{Thm}
\begin{proof}
Let $P\in\calVe(\calM)$ be the MC element with derived brackets $\lambda$. Define $m$ by applying the
$L_\infty$\ndash quasi\ndash isomorphism $U$ to $\epsilon P$:
\[
m = \sum_{n=0}^\infty \frac 1{n!} U_n(P,\dots,P).
\]
This is a MC element in $\epsilon\calDe(\calM)$ with the desired properties.
%
\end{proof}


Observe however that the problem of quantizing flat $P_\infty$\ndash manifolds in the flat sense
is on the other hand not solved. 
In fact, we are not able to conclude that
$P_0=0$ implies $m_0=0$ but only that
$m_0=O(\epsilon^2)$. This is the deformation quantization version
of a (potential) anomaly (see \cite{CF07}).

There are some cases when $m_0$ actually
happens to vanish (see \cite{CF04}). There are also cases \cite{CF07}
where $m_0$ can be killed by
shifting the $A_\infty$\ndash structure without changing its classical
limit. Namely, for $\gamma\in\epsilon A^1[[\epsilon]]$ 
define
\begin{align*}
\hat m_0 &= m_0 + m_1(\gamma) + \frac12 m_2(\gamma,\gamma) +\cdots,\\
\hat m_1(a) &= m_1(a) + m_2(a,\gamma) \pm m_2(\gamma,a) + \cdots.
\end{align*}
Then, $\forall\gamma$, 
$\mu+\hat m$ is again an $A_\infty$\ndash structure on $A[[\epsilon]]$
that induces the same $P_\infty$\ndash structure on $A$. One may
then look for a $\gamma$ such that $\hat m_0=0$. It is not difficult to prove
by induction, see \cite{CF07}, that a sufficient condition for this to happen is that
$\sfA^2=\{0\}$. So we have the following
\begin{Thm}\label{A2}
Every flat $P_\infty$\ndash manifold $\calM$ with $\sfA^2=\{0\}$
admits a deformation quantization in the flat sense.
\end{Thm}

\subsection{Quantization of cohomology}\label{s:qc}
Assume now that the flat $P_\infty$\ndash manifold $\calM$ has been quantized
in the flat $A_\infty$ sense. We denote by $\chi$ the graded commutative multiplication on the algebra of functions
$A$, by
$\lambda_i$ the multibrackets defining the
flat $P_\infty$\ndash structure on  $A$ and by
$\mu_i$ the multibrackets defining the  flat
$A_\infty$\ndash structure on $A[[\epsilon]]$. We will also denote by $\dd$ the differential $\lambda_1$ and
by $\sfA$ the $\dd$\ndash cohomology of $A$. 

The $\mu_1$\ndash cohomology 
is an associative algebra. However, in view of quantizing $\sfA^\bullet$, or at least $\sfA^0$,
this is not the algebra we are in general interested in 
since the projection
modulo $\epsilon$ is not a chain map.
We have instead to proceed as follows. Observe that by appropriately rescaling the
multibrackets, we get a new $A_\infty$\ndash structure; viz., for given integers $s_i$,
define
$\tau_i = \epsilon^{s_i}\mu_i$, 
where we assume $\mu_i$ to be divisible by $\epsilon^{-s_j}$ whenever $s_j<0$.
If moreover $s_i+s_j$ is a function of $i+j$, the multibrackets $\tau_i$ also define
an $A_\infty$\ndash algebra.

In particular, we may take $s_i=i-2$ since $\mu_1=O(\epsilon)$. Namely, we define the new 
$A_\infty$\ndash algebra structure
\[
\tau_i:=\epsilon^{i-2}\mu_i.
\]
Observe that
\[
\tau_1=\dd + O(\epsilon),\qquad
\tau_2=\chi + O(\epsilon),\qquad
\tau_i = O(\epsilon^{i-1}),\ i>2.
\]
By construction, the graded
skew\ndash symmetrization $\Tilde\tau_2$ of $\tau_2$ is divisible 
by $\epsilon$; so there is a unique operation $\psi$ with
$\Tilde\tau_2=\epsilon\psi$. The bracket $\lambda_2$ on $A$ is then $\psi$ modulo
$\epsilon$.

We will also denote by $\delta=\sum_{n=0}^\infty \epsilon \dd_n$, $\dd_0=\dd$, 
the differential $\tau_1$ and
by $\sfB$
the $\delta$\ndash cohomology 
of $A[[\epsilon]]$. Observe that the differential $\delta$ is a deformation 
of the differential $\dd$. 

Observe that by construction $\tau_2$, $\Tilde\tau_2$ and $\psi$
are chain maps and we will denote by $[\tau_2]$, $[\Tilde\tau_2]$ and $[\psi]$
the induced operations in cohomology. The operation $[\tau_2]$ is an
associative product which we will also
denote by $\sfstar$
(and by $\star$ its restriction to $\sfB^0$).

Whenever needed we denote $\delta$\ndash cohomology classes by $[\ ]_\delta$ 
and $\dd$\ndash cohomology classes by $[\ ]_\dd$. 
With these notations we have
\[
[a]_\delta\sfstar [b]_\delta = [\tau_2]([a]_\delta, [b]_\delta)
= [\tau_2(a,b)]_\delta,
\]
where $[\tau_2]$ is the map induced by $\tau_2$ in cohomology, while
$a$ and $b$ are representatives of the classes $[a]_\delta$ and $[b]_\delta$.
By construction $[\Tilde\tau_2]$ defines
the graded $\sfstar$\ndash commutator $\sfLie{\ }{\ }$ on $\sfB$:
\[
\sfLie{[a]_\delta}{[b]_\delta} = [\Tilde\tau_2]([a]_\delta, [b]_\delta)
= [\Tilde\tau_2(a,b)]_\delta.
\]
We will denote simply by $\Lie{\ }{\ }$ its restriction to $\sfB^0$.
Observe that multiplication by $\epsilon$ commutes with $\delta$. 
So
$\sfB$ is also an
$\bbR[[\epsilon]]$\ndash module with $\epsilon[a]_\delta=[\epsilon a]_\delta$.
Moreover, we have
\[
\sfLie{[a]_\delta}{[b]_\delta} =\epsilon[\psi]([a]_\delta, [b]_\delta).
\]
If $\sfB$ ($\sfB^0$) is $\epsilon$\ndash torsion 
free---i.e., multiplication by $\epsilon$ 
is injective---, then $[\psi]$ is the unique operation with the above property.
Consider now the $\bbR$\ndash linear projection
\[
\varpi\colon\begin{array}[t]{ccc}
A[[\epsilon]] &\to &A\\
\sum_{n=0}^\infty \epsilon^n a_n &\mapsto &a_0.
\end{array}
\]
Observe that 
$\varpi\colon (A[[\epsilon]],\delta)\to(A,\dd)$ is a chain map.
We will denote by $[\varpi]\colon\sfB\to\sfA$ the induced map in cohomology.
Since $\varpi$ is an algebra homomorphism, so is
$[\varpi]$.
The relation with the Poisson structure on $\sfA$ is clarified by the now
obvious formula
\[
[\varpi]([\psi]([a]_\delta,[b]_\delta)) =
\sfPoiss{[\varpi]([a]_\delta)}{[\varpi]([b]_\delta)}, \qquad
\forall [a]_\delta,[b]_\delta\in \sfB
\]
which immediately implies the
\begin{Prop}\label{dequant} 
The algebra homomorphism $[\varpi]$ has the following additional properties:
\begin{enumerate}
\item The image $\sfC$ ($\sfC^0$) of $[\varpi]$ is a Poisson subalgebra of
$\sfA$ ($\sfA^0$). 
\item If $\sfB$ is $\epsilon$\ndash torsion free, then 
\[
[\varpi]\left(\frac{\sfLie{[a]_\delta}{[b]_\delta}}\epsilon\right) =
\sfPoiss{[\varpi]([a]_\delta)}{[\varpi]([b]_\delta)}, \qquad
\forall [a]_\delta,[b]_\delta\in \sfB.
\]
\item If $\sfB^0$ is $\epsilon$\ndash torsion free, then the above formula
holds for all $[a]_\delta,[b]_\delta\in \sfB^0$.
\end{enumerate}
\end{Prop}
Since obviously
$\Ker\varpi=\epsilon A[[\epsilon]]$, we have
the  short exact sequence
\[
0\rightarrow A[[\epsilon]] \xrightarrow\epsilon A[[\epsilon]]
\xrightarrow\varpi A\rightarrow 0
\]
which induces the long exact sequence 
\[
\cdots\rightarrow  \sfA^{i-1} \xrightarrow\partial
\sfB^i \xrightarrow\epsilon \sfB^i \xrightarrow{[\varpi]} \sfA^i
\xrightarrow\partial \sfB^{i+1} \rightarrow\cdots
\]
in cohomology.
Immediately we then get the
\begin{Prop}\label{equi}
The following statements are equivalent:
\begin{enumerate}
\item $[\varpi]$ is surjective;
\item $\sfB$ is $\epsilon$\ndash torsion free;
\item $\partial$ is trivial.
\end{enumerate}
\end{Prop}

To continue our study of the problem, we now need 
Lemma~A.1 of \cite{CFT}, which we state in a slightly modified version:
\begin{Lem}\label{lemCFT}
Let $\Bbbk$ be a field and $\sfM$ a $\Bbbk[[\epsilon]]$\ndash module endowed
with the $\Bbbk[[\epsilon]]$\ndash adic topology. 
Then 
$\sfM\simeq_\Bbbk\sfM_0[[\epsilon]]$ for some $\Bbbk$\ndash vector space $\sfM_0$
if{f} $\sfM$ is Hausdorff, complete and
$\epsilon$\ndash torsion free.
Moreover, $\sfM_0\simeq_\Bbbk \sfM/\epsilon\sfM$.
\end{Lem}
We finally have the
\begin{Thm}\label{fawkes}
If $\sfB$ is Hausdorff and complete in the $\epsilon$\ndash adic topology, then
$\sfB$ is a deformation quantization of $\sfA$ if{f} any (and so all) of the
statements in Proposition~\ref{equi} holds.
\end{Thm}
\begin{proof}
If $\sfB$ is a deformation quantization of $\sfA$, then
in particular
$\sfB\simeq_\bbR\sfA[[\epsilon]]$. So by Lemma~\ref{lemCFT} $\sfB$
is $\epsilon$\ndash torsion free.
On the other hand, the statements in  Proposition~\ref{equi} imply that
$\sfA\simeq_\bbR\sfB/\Ker[\varpi]=\sfB/\epsilon\sfB$; so
$\sfB\simeq_\bbR\sfA[[\epsilon]]$ by Lemma~\ref{lemCFT}. 
Statement (2) in Proposition~\ref{dequant} completes the proof.
\end{proof}

We are not able to show that $\Ima\delta$ is closed in the 
$\epsilon$\ndash adic topology,
so we must put the extra condition
in the Theorem. We wish however to make the following\footnote{Observe that the Conjecture is general
false if we drop the conditions that $A$ is the algebra of functions
on a graded manifold and that the components $\dd_n$ of 
$\delta$ are differential operators.} 
\begin{Conj}\label{conj}
The $\bbR[[\epsilon]]$\ndash module
$\sfB$ is Hausdorff and complete in the $\epsilon$\ndash adic topology.
\end{Conj}


In general 
we do not expect $\sfB$ to be a deformation quantization of $\sfA$.
However, it is often enough to have $\sfB^0$ as a deformation quantization
of $\sfA^0$. We end this Section by exploring some sufficient conditions for
this to happen. 
\begin{Lem}\label{tfree}
If $\sfB^i$ is Hausdorff and complete, then
it is isomorphic to $\sfC^i[[\epsilon]]$
as an $\bbR$\ndash vector space if{f}
it is $\epsilon$\ndash torsion free.
\end{Lem}
\begin{proof}
The ``only if'' implication is obvious. As for the ``if'' part,
by Lemma~\ref{lemCFT} we conclude that  $\sfB^i$
is isomorphic to $(\sfB^i/\epsilon\sfB^i)[[\epsilon]]$.
On the other hand the long exact sequence gives 
$0\rightarrow
\sfB^i \xrightarrow\epsilon \sfB^i \xrightarrow{[\varpi]} \sfC^i
\rightarrow0$
which completes the proof.
\end{proof}
By statement (3) in Proposition~\ref{dequant}, we then get
\begin{Cor}\label{cordq}
If $\sfB^0$ is Hausdorff, complete and
$\epsilon$\ndash torsion free, then it is a deformation quantization
of $\sfC^0$.
\end{Cor}
The long exact sequence yields the
\begin{Lem}\label{Am1}
If $\sfA^{i-1}=\{0\}$, then $\sfB^i$ is $\epsilon$\ndash torsion free.
\end{Lem}
Thus, $A^{-1}=\{0\}$ is a sufficient condition for $\sfB^0$ to be a deformation
quantization of $\sfC^0$ provided it is Hausdorff and complete.
Finally, we have the
\begin{Lem}\label{A1}
If $\sfA^{i+1}=\{0\}$, 
then $[\varpi]^i\colon\sfB^i\to\sfA^i$ is surjective.
\end{Lem}
\begin{proof}
This is a standard proof in cohomological perturbation theory.
Let $[a_0]_\dd\in\sfA^i$. Choose one of its representatives $a_0\in A^i$.
We look for a $\delta$\ndash closed 
$a=\sum_{n=0}^\infty\epsilon^n a_n\in A^i[[\epsilon]]$.
The equations we have to solve have the form
\begin{equation}\label{dda}
\dd a_n = -\sum_{\substack{r+s=n+1\\ r>0}}
\dd_r a_s,
\end{equation}
and we may solve them by induction. Namely:
\begin{enumerate}
\item For $n=0$, the equation is just $\dd a_0=0$ which is satisfied by assumption.
\item Assume now that all equations for $a_k$, $k<n$,
have been solved. This implies that that the r.h.s.\ of \eqref{dda}
is $\dd$\ndash closed. Since we assume
that $\sfA^{i+1}=\{0\}$, we may then find $a_n$ satisfying the equation.
\end{enumerate}
\end{proof}

\begin{Rem}
Observe that, if we knew that $B^{i+1}$ were Hausdorff (e.g., if 
Conjecture~\ref{conj} were true), then we could derive Lemma~\ref{A1}
directly from the long exact sequence (this is just a  variant of Nakayama's Lemma).
In fact, $A^{i-1}=\{0\}$
implies that multiplication by $\epsilon$ on 
$\sfB^{i+1}$ is surjective. So every element $a\in\sfB^{i+1}$ may be written
as $a=\epsilon a_1$. Continuing this process, we get a sequence $a_n$ with
$a=\epsilon^n a_n$. The sequence $\epsilon^n a_n$ converges to zero, 
so $a=0$ since
$\sfB^{i+1}$ is Hausdorff. Thus, $\sfB^{i+1}=\{0\}$ and by the long exact 
sequence again
we get the result.
\end{Rem}



Putting together Theorem~\ref{A2}, Corollary~\ref{cordq},
Lemmata~\ref{Am1}  and~\ref{A1},
we get the following
\begin{Thm}\label{auror}
If 
$\sfA^2=\sfA^1=\sfA^{-1}=\{0\}$, 
then $\sfB^0$ is 
a deformation quantization of $\sfA^0$ provided it is Hausdorff and complete.
\end{Thm}

If Conjecture~\ref{conj} holds, we get the sufficient conditions in a much nicer way
which makes reference only to the $\dd$\ndash cohomology:
\begin{Thmt}\label{squib}
A sufficient condition for a deformation quantization of $\sfA^0$ to exist is
$\sfA^2=\sfA^1=\sfA^{-1}=\{0\}$.
\end{Thmt}

Observe that the sufficient condition is by no means necessary.
For example, in the case $\lambda_i=0$ $\forall i$, we have $\sfA^\bullet=A^\bullet$.
On the other hand, $\delta=0$ and $\sfB^\bullet=A^\bullet[[\epsilon]]$ 
is a deformation quantization
of $\sfA^\bullet$.


\subsection{Deformation quantization of coisotropic submanifolds}\label{s:dqcs}
We now return to the case when
$C$ is a coisotropic submanifold of $M$ (which is Poisson at least
in a neighborhood of $C$). We take our graded manifold $\calM$ to be $N^*[1]C$ with MC element
the restriction of $\pi$ to $\calV(NC)=\calV(N^*[1]C)$. 
Now $\sfA^{-1}=\{0\}$ since we do not have negative degrees, while 
$\sfA^0$ is the Poisson algebra
of functions on $C$ that are invariant under the canonical distribution on $C$.
If we denote by $\underline C$ the leaf space, we also write
$C^\infty(\underline C)=\sfA^0$.
Observe that in general $\underline C$ is not a manifold, so the above is just a definition.
Using the notations of Section~\ref{s:qc}, we have $\sfB^0=\Ker\delta$. So $\sfB^0$
is Hausdorff and complete, and by Theorem~\ref{auror} we have the
\begin{Thm}[\cite{CF07}]\label{t:noanomaly}
If the first and second Lie algebroid cohomology of $N^*C$ vanish,
the zeroth cohomology $C^\infty(\underline C)$ has a deformation quantization.
\end{Thm}
The second Lie algebroid cohomology is also the space where the usual BRST/BFV anomaly lives; see
\cite{B04} for this in the context of deformation quantization. 
The advantage of the present approach
is that one has formulae for $m_0$, so one can in principle check whether it may be canceled.

Observe that the conditions of Theorem~\ref{t:noanomaly} are very strong and by no means necessary. At the moment we actually do not know
a single example where the anomaly shows up, at least locally. 
There are examples where the Lie algebroid cohomology class of $m_0$ is nontrivial, see \cite{W},
but in these examples $C$ is a point, so there is no problem in quantizing its reduced space (the problem arises however in the many-brane setting
of Section~\ref{s:MB}).

As $C^\infty(\underline C)$ might be rather poor (e.g., just constant functions),
it might be interesting to quantize the whole Lie algebroid cohomology. Assuming that the 
$A_\infty$\ndash structure is flat,
by Theorem~\ref{fawkes}
the deformed cohomology $\sfB$ yields a deformation quantization of $\sfA$
if{f} it is Hausdorff, complete and 
$\epsilon$\ndash torsion free.

\subsection{Poisson submanifolds}\label{s:Psub}
A particular case of a coisotropic submanifold is a Poisson submanifold, i.e.,  a submanifold $C$ of $(M,\pi)$ such that the inclusion map is a Poisson map.
Equivalently, a Poisson submanifold is a coisotropic submanifold with trivial characteristic distribution.
In this case, $\und C =C$, so the quantization of the reduction is not a problem. The interesting question is whether one can deform the pullback of the inclusion
map to a morphism of associative algebras. To approach this problem, we associate to $C$ a different graded manifold as in the coisotropic case.

If the submanifold $C$ is determined by constraints, i.e., $C=\Phi^{-1}(0)$ with $\Phi\colon M\to V$ a given map to a vector space $V$, 
we set $\calM:=M\times V[-1]$ and reinterpret $\Phi$
as an element of $C^\infty(M)\otimes V$ and so
as a vector field $X$ of degree $1$ on $\calM$. If we introduce coordinates $\{\mu^\alpha\}_{\alpha=1,\dots,k:=\dim V}$ of degree $-1$, we have 
$X=\Phi^\alpha(x)\frac\de{\de\mu^\alpha}$,
where the $\Phi^\alpha$s are the components of $\Phi$ w.r.t.\ the same basis. Observe that $\Lie XX=0$ and that $\Lie\pi X$ vanishes on $C$.
Moreover, the cohomology of $C^\infty(\calM)$ w.r.t\ the differential $\delta:=\Lie X{\ }$ is $C^\infty(C)$.
If $C$ is not determined by constraints, we take $\calM:=N[-1]C$ and $X$ the vector field of degree $1$ corresponding to the zero section.
The crucial observation now is that $A:=C^\infty(\calM)=\bigoplus_{j\le0} A^j$ is a nonpositively graded commutative algebra. As a consequence,
$\sfA^i=0$, $i>0$,  with the differential $\lambda_1=\delta$.
Observe moreover that $A^0$ is just $C^\infty(M)$ and that $\sfA^0=A^0/I=C^\infty(C)$, with $I=\delta A^{-1}$ the vanishing ideal of $C$.

\subsubsection{Casimir functions}
A particularly simple case is when the components of $\Phi$ are Casimir functions: viz., $\Poiss{\Phi^\alpha}f=0$ $\forall\alpha$ $\forall f\in C^\infty(M)$.
Equivalently $\Lie\pi{\Phi^\alpha}=0$ $\forall\alpha$, i.e., $\Lie\pi X=0$. Thus $\Hat\pi:=\pi+X$ is a MC element. The induced $P_\infty$\ndash structure
has $\lambda_1=\delta$. Quantization has the following easy to check remarkable properties: $i)$ $(A^0[[\epsilon]],\tau_2)$ is Kontsevich's deformation 
quantization $A_M$ of $C^\infty(M)$;
$ii)$ $\tau_2(\mu^\alpha,f)=\tau_2(f,\mu^\alpha)=\mu^\alpha f$, $\forall\alpha$. As a consequence $\sfI:=\tau_1(A^{-1}[[\epsilon]])$ is the two-sided ideal
generated by $\{\tau_1(\mu^\alpha)\}$. A further easy computation, using Stokes' theorem in the explicit 
expression of the coefficients by Kontsevich's graphs, shows that $\tau_1(\mu^\alpha)=D(\Phi^\alpha)$, where $D$ is the
Duflo--Kirillov--Kontsevich map
\[
D\colon\begin{array}[t]{ccc}
C^\infty(M) &\to & C^\infty(M)[[\epsilon]],\\
f &\mapsto & \sum_{n=0}^\infty \frac{\epsilon^n}{n!}\,U_{n+1}(f,\pi,\dots,\pi).
\end{array}
\]
So we get a deformation quantization of $C$ as $\sfB^0=A_M/\sfI$, and the projection $A_M\to\sfB^0$ is a quantization
of the inclusion map.
This proves a conjecture in \cite{CR} (see also \cite[Sect.~5]{C06} for a previous sketch of this proof).
\begin{figure}[ht] 
\begin{center}
\resizebox{8 cm}{!}
{\includegraphics{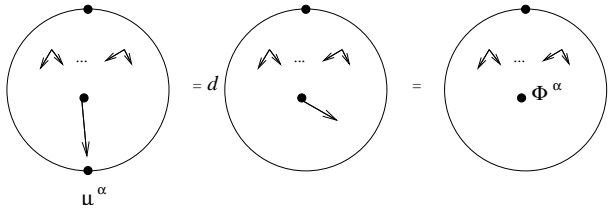}}
\caption{$\tau_1(\mu^\alpha)=D(\Phi^\alpha)$}  
\end{center}
\end{figure}

\subsubsection{The general case}
In general $\pi+X$ is not a MC element. However, the fact that $\Lie\pi X$ vanishes on $C$, makes $\pi+X$ a MC element up to $\delta$\ndash exact terms.
As a consequence, one may use cohomological perturbation theory as in \cite{LS} and find a MC element $\hat\pi$ of the form $\pi+X+$corrections in the 
ideal of multivector fields of total degree $2$ generated by $\{\mu^\alpha\}$.
Observe that the only vector field in $\hat\pi$ is $X$, so $\lambda_1=\delta$. The vanishing of the cohomologies in positive degrees implies that we have
an $A_\infty$\ndash structure on $A[[\epsilon]]$ and a surjective map $\sfB^0\to A^0/I$. Moreover, the restriction of $\tau_2$ in degree $0$ makes
$A^0[[\epsilon]]$ into an algebra with an algebra morphism to $\sfB^0=A^0[[\epsilon]]/\sfI$, $\sfI:=\tau_1(A^{-1}[[\epsilon]])$.
The two problems of this constructions are the following: $i)$  $(A^0[[\epsilon]],\tau_2)$ might not be associative;
$ii)$ $\sfB^0$ might not be isomorphic to $A^0/I[[\epsilon]]$.

It may be shown that $(A^0[[\epsilon]],\tau_2)$ is associative if{f} $\hat\pi$ has the form $\pi+X+\pi'$ with $\pi'$ 
a \emph{bivector}\/ field. 
If in addition one chooses the constraints to be linear (e.g.,
if one works with $\calM=N[-1]C$), then one can also  easily see that
$\sfB^0$ turns out to be a deformation quantization of $C^\infty(C)$.
However, finding a $\pi'$ as above is a highly nontrivial problem, and it is not clear under which conditions a solution may exist.

A very simple case is when $C$ consists of a point $x$ (a zero of the Poisson structure). A quantization of the inclusion map can then be reinterpreted as 
a character (i.e., an algebra morphism to the ground ring) of the deformed algebra (deforming the evaluation at $x$).
Even in such a simple situation, the existence of a $\pi'$ is guaranteed only for $\dim M=2$ \cite{S2}.
In higher dimensions, it is an open problem.



\section{Many branes}\label{s:MB}
We now turn to the case when more than one coisotropic submanifold is chosen as a boundary condition.
Namely, as in Section~\ref{s:exp},
we take $\Sigma$ to be a disk. However, we now subdivide the boundary into closed, cyclically ordered intervals $I_1,\dots,I_n$ with exactly one intersection point
between subsequent intervals. To the interval $I_i$ we associate boundary conditions corresponding to a coisotropic submanifold $C_i$.

Assuming clean intersections,
in \cite{CF04} it is shown that the case $i=2$ leads to the construction of a bimodule for the deformation quantizations of $\und C_1$ and $\und C_2$ if no anomaly
appears.

In this Section we assume that no anomalies show up. If the submanifolds are pairwise transverse this amounts to asking that for each of them
the anomaly discussed in subsection~\ref{s:dqcs} 
vanishes. The vanishing of the second Lie algebroid cohomology for each submanifold is a sufficient condition
by Theorem~\ref{t:noanomaly}. Notice however that this is only a sufficient condition by no means necessary.
A very simple example is when $M$ is symplectic and $C_1$ is Lagrangian. In this case $\und C_1$ is a point (assume $C_1$ to be connected), so there
is no problem in quantizing it. The bimodule structure associated to $C_2=M$ can always be found (upon choosing the star product appropriately).
On the other hand the relevant Lie algebroid cohomology is the de~Rham cohomology of $C_1$. Another example
is the linear case $C_1=\frh^0\subset \frg^*$, where $\frh$ is a Lie subalgebra of $\frg$. It is shown in \cite{CF04} that there is no anomaly 
in this case even if the cohomology may be very complicated. Observe however that an example of nonvanishing anomaly has recently been found in \cite{W}.

\begin{figure}[ht] 
\begin{center}
{\includegraphics{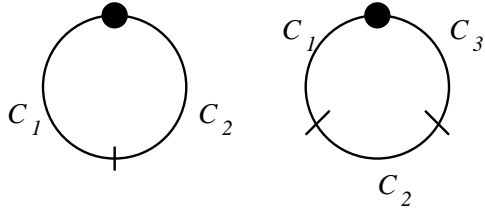}}
\caption{Two and three branes}  \label{f:bimod}
\end{center}
\end{figure}

We also use the notation $A_{\und C}$ for the PSM deformation quantization of the Poisson algebra
$C^\infty(\und C)$, where $C$ is a coisotropic submanifold of $M$. Whenever a distribution is trivial, we omit underlining the submanifold.
Observe that $A_M$ is Kontsevich's deformation quantization.
With $B_{\und{C_1\cap C_2}}$ we then denote the $A_{\und{C_1}}$\ndash module\ndash $A_{\und{C_2}}$ corresponding to the picture on the left in fig.~\ref{f:bimod}
(where the big black dot denotes the point at infinity which has to be sent to the point in the intersection where evaluation takes place).

The situation with three coisotropic branes $C_1$, $C_2$ and $C_3$---see the picture on the right in fig.~\ref{f:bimod}---leads to a morphism of bimodules 
$B_{\und{C_1\cap C_2}}\otimes_{A_{\und{C_2}}}B_{\und{C_2\cap C_3}}\to B_{\und{C_1\cap C_3}}$ in a neighborhood of a triple intersection. The construction
of the corresponding Kontsevich graphs has been analyzed in \cite{CT}. In subsection~\ref{s:fq} we discuss a special case.

One may also consider many branes and, instead of just invariant functions, sections of the corresponding complexes of exterior algebras of normal bundles.
One may hope that this would lead to $A_\infty$\ndash bimodules and morphisms thereof, but problems seem to arise with four and more branes.

In the rest of this Section we will concentrate on the special case when one of the branes is the whole manifold; we will consider only functions and
no more than three branes.

\subsection{Quantization of morphisms}
Let $\phi$ be a Poisson map $M\to N$. Then its graph $\Graph\phi$ is a coisotropic submanifold of  $\overline M\times N$,
where $\overline M$ denotes $M$ with opposite Poisson structure. We select the two branes $C_1=\Graph\phi$ and $C_2=M\times N$
in $\overline M\times N$.
Consider the  right
$A_{\overline M\times N}$\ndash module structure on $B_{\Graph\phi}$, forgetting about its left $A_{\und{\Graph\phi}}$\ndash module structure.
Since $A_{\overline M}\otimes A_N$ is a subalgebra of $A_{\overline M\times N}$ and  Kontsevich's star product has the property
$A_{\overline M}=A_M^\mathrm{opp}$, we may regard $B_{\Graph\phi}$ as an $A_M$\ndash bimodule\ndash$A_N$. This bimodule has a distinguished
element, the constant function $1$, and the map $A_M\to B_{\Graph\phi}$, $f\mapsto f\cdot 1$ is an isomorphism of $\bbR[[\epsilon]]$\ndash modules
(since it is a deformation of the pullback of the diffeomorphism $p_1\colon\Graph\Phi\to M$ defined by projection). Thus, to every $f\in A_N$
we associate a unique element $\Hat\phi(f)$ of $A_M$ by the equation
\[
1\cdot f = \Hat\phi(f)\cdot 1.
\]
Observe that $\Hat\phi$ is a deformation of the pullback $\phi^*$. It is not difficult to see that it is a morphism of associative algebras.
This way, with the assumption that there is no anomaly, we have found a quantization procedure for Poisson maps.

Observe that one can exchange the role of $C_1$ and $C_2$ using the symmetry corresponding to reflecting the disk through the line joining the two 
intersections of the intervals $I_1$ and $I_2$ ($\bullet$ and $|$ in fig.~\ref{f:bimod}).
In terms of Kontsevich graphs one gets the same formulae upon changing the sign of the Poisson structure;
viz., we get the same result as before if we take $M\times\overline N$ as the ambient Poisson manifold. 

\subsection{Compositions}\label{s:fq}
Now suppose we have two Poisson maps $\phi\colon M\to N$ and $\psi\colon N\to K$. With no anomalies, we may quantize $\phi$, $\psi$ and $\psi\circ\phi$.
The next problem we wish to address concerns the relation between $\widehat{\psi\circ\phi}$ and $\Hat\phi\circ\Hat\psi$.

To deal with it we consider the three branes $C_1=\Graph\phi\times K$, $C_2=M\times N\times K$ and $C_3=M\times\Graph\psi$ in $\overline M\times N\times\overline K$
as in fig.\ref{f:bimod}. Let $p_1\colon M\times N\to M$ and $\tilde p_2\colon N\times K\to K$ be the canonical projections.
To $f\in C^\infty(M)[[\epsilon]]$ and $g\in C^\infty(K)[[\epsilon]]$, we associate the element $D(f,g)\in C^\infty(\Graph\psi\circ\phi)[[\epsilon]]$
obtained by putting $p_1^*f\otimes1$ and $1\otimes\tilde p_2^*g$ at the two intersection points (here $1$ is the constant function on $C^\infty(N)$).
It is possible to check that $D$ defines a morphism of  $A_M$\ndash bimodules\ndash$A_K$ 
$B_{\Graph\phi}\otimes_{A_N}\otimes B_{\Graph\psi}\to B_{\Graph\psi\circ\phi}$. 
So it is enough to compute $\sigma_{\phi,\psi}:=D(1,1)=1+O(\epsilon^2)\in B_{\Graph\psi\circ\phi}$.
Observe again that there is a unique element $\Hat\sigma_{\phi,\psi}\in A_M$ such that $\Hat\sigma_{\phi,\psi}\cdot 1=\sigma_{\phi,\psi}$. Moreover, since $\Hat\sigma_{\phi,\psi}$ is of the form
$1+O(\epsilon^2)$, it is invertible. Now for $h\in A_K$ we have the identities
\begin{multline*}
\sigma_{\phi,\psi}\cdot h 
 = D(1,1\cdot h) = D(1,\Hat\psi(h)\cdot 1)=\\
=D(1\cdot\Hat\psi(h),1)=D(\Hat\phi(\Hat\psi(h))\cdot1,1)=
\Hat\phi\circ\Hat\psi\cdot\sigma_{\phi,\psi}.
\end{multline*}
Using $1\cdot h=\widehat{\psi\circ\phi}(h)\cdot 1$ and the definition of
$\hat\sigma_{\phi,\psi}$, we get
\begin{gather*}
\sigma_{\phi,\psi}\cdot h = \Hat\sigma_{\phi,\psi}\cdot 1\cdot h = \Hat\sigma_{\phi,\psi}\cdot(\widehat{\psi\circ\phi}(h)\cdot 1)=
(\Hat\sigma_{\phi,\psi}\star_M \widehat{\psi\circ\phi}(h))\cdot 1,\\
\Hat\phi\circ\Hat\psi(h)\cdot\sigma_{\phi,\psi}=\Hat\phi\circ\Hat\psi(h)\cdot(\Hat\sigma_{\phi,\psi}\cdot 1)=
(\Hat\phi\circ\Hat\psi(h)\star_M\Hat\sigma_{\phi,\psi})\cdot 1,
\end{gather*}
where $\star_M$ denotes the star product on $M$. Finally,
\[
\widehat{\psi\circ\phi}(h)=\Hat\sigma_{\phi,\psi}^{-1}\star_M (\Hat\phi\circ\Hat\psi(h))\star_M\Hat\sigma_{\phi,\psi},
\]
where $\hat\sigma_{\phi,\psi}^{-1}$ is the inverse of $\hat\sigma_{\phi,\psi}$ w.r.t.\ $\star_M$.

Thus, upon conjugation (by an element which depends on the given Poisson maps), we see that the composition of quantizations is the quantization of compositions.

\subsection{Quantum groups}
We now want to apply the results of the last subsection to the case of a Poisson--Lie group $G$; i.e., a Poisson manifold and a Lie group such that the
product $m\colon G\times G\to G$ and the inverse ${}^{-1}\colon \overline G\to G$ are
Poisson maps. By Kontsevich we have an associative algebra $A_G$. On the other hand, if no anomaly is present,
we get a morphism of algebras $\Delta:=\Hat m\colon A_G\to A_G\hat\otimes A_G:=A_{G\times G}$.
The associativity equation $m\circ (m\otimes\mathit{id})=m\circ(\mathit{id}\otimes m)$ then yields
\[
(\Delta\otimes\mathit{id})\circ\Delta = \Phi^{-1}\star_{G^3}((\mathit{id}\otimes\Delta)\circ\Delta)\star_{G^3}\Phi,
\]
with $\Phi=\hat\sigma_{\mathit{id}\otimes m,m}^{-1}\star_{G^3}\hat\sigma_{m\otimes\mathit{id},m}\in A_{G^3}$.
In some lucky cases, $\Phi$ might turn out to be $1$ (or at least central). In this case we would have quantized the Poisson--Lie group as a bialgebra. Otherwise, one may
hope to get the right properties for the  ``associator'' $\Phi$ in order to get a bialgebra out of it.

Finally, observe that we may also quantize the inverse and get a candidate for the antipode. If relations are preserved, we should get a Hopf algebra structure
on $A_G$, i.e., the corresponding quantum group. 

Just assuming that anomalies are absent is not enough to get a Hopf algebra as relations are preserved only up to conjugation by certain elements.
If one cannot get rid of them, the resulting structure is probably that of a hopfish algebra as defined in \cite{BW}.






\bibliographystyle{amsalpha}

\end{document}